\documentclass[12pt]{article}
\usepackage{stmaryrd}
\usepackage{mathrsfs}
\usepackage{amsmath}
\usepackage{amsmath,amsthm,amssymb}
\usepackage{amsfonts}
\usepackage{latexsym}
\date{}

\begin{document}
\title{Tilting modules over duplicated algebras}
\author{{\small  Guopeng Wang, Shunhua Zhang}\\
{\small  School of Mathematics,\ Shandong University,\ Jinan 250100,
P. R. China }}

\pagenumbering{arabic}

\maketitle
\begin{center}
 \begin{minipage}{120mm}
   \small\rm
   {\bf  Abstract}\ \ Let $A$ be a finite dimensional hereditary
algebra over a field $k$ and $A^{(1)}$ the duplicated algebra of
$A$. We first show that the global dimension of endomorphism ring of
tilting modules of $A^{(1)}$ is at most 3. Then we investigate
embedding tilting quiver $\mathscr{K}(A)$ of $A$ into tilting quiver
$\mathscr{K}(A^{(1)})$ of $A^{(1)}$.  As applications, we give new
proofs for some results of D.Happel and L.Unger, and prove that
every connected component in $\mathscr{K}({A})$ has finite
non-saturated points if $A$ is tame type, which gives a partially
positive answer to the conjecture of D.Happel and L.Unger in [10].
Finally, we also prove that the number of arrows in
$\mathscr{K}({A})$ is a constant which does not depend on the
orientation of $Q$ if $Q$ is Dynkin type.

\end{minipage}
\end{center}

\footnote{ {\it Email addresses}: \
guopengwang@mail.sdu.edu.cn(P.Wang),
 \ shzhang@sdu.edu.cn(S.Zhang)}

\vskip0.2in

\section {Introduction}

Tilting theory usually has two aspects. One is the external aspect,
which is used to compare $\Lambda$-mod to ${\rm End}_{\Lambda}
T$-mod for a tilting $\Lambda$ module $T$. The other internal
aspect, which is to study tilting modules for a fixed algebra
$\Lambda$ and to try to gather information about $\Lambda$-mod, see
[6, 7, 8, 9, 10] for more details. Recently, tilting theory has
remarkable development in cluster categories, which was introduced
in [5].  Now, cluster categories become a successful model for
acyclic cluster algebras, this new discovery has rapidly promoted
research on this direction.

\vskip0.2in

According to [1], we know that tilting modules of duplicated algebra
$A^{(1)}$ of hereditary algebra $A$ have strong relationship with
cluster tilting objects in cluster category ${\cal C}_A$. For
example,  there is a one-to-one correspondence between basic tilting
$A^{(1)}$-modules with projective dimension at most one and basic
cluster tilting objects in ${\cal C}_A$.

\vskip0.2in

It is well known that the tilting quiver $\mathscr{K}(A)$ of a
hereditary algebra $A$ usually is not connected. For example, the
tilting quiver of  Kronecker algebra consists of two connected
components. However, according to [14], we know that the tilting
quiver $\mathscr{K}(A^{(1)})$ of $A$ is connected. This motivates
further investigation on the structure of tilting modules of
hereditary algebras and of tilting modules of duplicated algebras
with projective dimension at most one.

\vskip0.2in

In this paper, we focus on the structure properties of tilting
modules with projective dimension at most one for duplicated algebra
$A^{(1)}$, and prove that the global dimension of endomorphism ring
of this kinds tilting modules is at most 3 (see Theorem 3.1 in
Section 3). In Section 4, we are interested in the relationship
between the tilting quivers  $\mathscr{K}(A)$ and
$\mathscr{K}(A^{(1)})$, and prove some embedding theorems.

\vskip0.2in

In Section 5, we give new proofs for some results of D.Happel and
L.Unger by using embedding theorem, and prove that every connected
component in $\mathscr{K}({A})$ has finite non-saturated points if
$A$ is tame type, which gives a partially positive answer to the
conjecture of D.Happel and L.Unger in [10]. We also prove that the
number of arrows in $\mathscr{K}({A})$ is a constant which does not
depend on the orientation of $Q$ if $Q$ is Dynkin type.  We fix
notations and recall some facts needed for our later use in Section
2.

\vskip0.2in

\section {preliminaries}

\vskip 0.2in

  Let $\Lambda$ be a finite dimensional
$k$-algebra over a field $k$ and $\Lambda$-mod be the category of
all finitely generated left $\Lambda$-modules. We denote by
$\Lambda$-ind the full subcategory of $\Lambda$-mod consisting of
indecomposable $\Lambda$ modules, and denote by ${\rm pd}_\Lambda X$
the projective dimension of an $\Lambda$ module $X$ and by ${\rm
gl.dim}\ \Lambda$ the global dimension of $\Lambda$. Let $D={\rm
Hom}_k(-,\ k)$ be the standard duality between $\Lambda$-mod and
$\Lambda^{op}$-mod, and $\tau_\Lambda$ be the Auslander-Reiten
translation of $\Lambda$. The Auslander-Reiten quiver of $\Lambda$
is denoted by $\Gamma_\Lambda$.

\vskip 0.2in

Given any module $M\in \Lambda $-mod, we denote by $M^{\bot}$ the
subcategory of $\Lambda$-mod with objects $X\in \Lambda$-mod
satisfying ${\rm Ext}_\Lambda^{i}(M,X)=0$ for all $i\geq 1 $  and by
$^{\perp}M$ the subcategory of $\Lambda$-mod with objects $X\in
\Lambda$-mod satisfying ${\rm Ext}_\Lambda^{i}(X,M)=0$ for $ i\geq
1$. We denote by $\Omega_{\Lambda}^{i}M$ and
$\Omega_{\Lambda}^{-i}M$ the i-th syzygy and cosyzygy  of $M$
respectively, and denote by gen $M$ the subcategory of $\Lambda$-mod
whose objects are generated by $M$. We may decompose $M$ as
$M\cong\oplus_{i=1}^{m}M_{i}^{d_{i}}$, where each $M_{i}$ is
indecomposable, $d_{i}>0$ for each $i$, and $M_{i}$ is not
isomorphic to $M_{j}$ if $i\neq j$. The module $M$ is called basic
if $d_{i}=1 $ for any $i$. The number of non-isomorphic
indecomposable modules occurring in the direct sum decomposition
above is uniquely determined and it is denoted by $\delta(M)$.  The
full subcategory having as objects the direct sums of indecomposable
summands of $M$ is denoted by ${\rm add} \ M$.

\vskip 0.2in

A module $T\in \Lambda$-mod is called a tilting module if
the following conditions are satisfied:\\
{(1)} ${\rm pd}_\Lambda T \leq 1$;\\
{(2)} ${\rm Ext}_\Lambda^{1}(T,T)=0 $;\\
{(3)} There is  an exact sequence
 $0\longrightarrow \Lambda \longrightarrow T_{0}\longrightarrow T_{1}
\longrightarrow 0$ with $T_{i}\in {\rm add}\ T $ for $0\leq i\leq 1$.

\vskip 0.2in

  An $\Lambda$ module $M$  satisfying the conditions $(1)$ and $(2)$ of
the definition above is called a partial tilting module and if
moreover $\delta(M)=\delta(\Lambda)-1 $, then  $M$ is called an
almost complete tilting module. Let $M$ be a partial tilting module
and $X$ be an $\Lambda$-module such that $M\oplus X$ is a tilting
module and ${\rm add}M\cap {\rm add}X=0$. Then X is called a
complement to $M$.

\vskip 0.2in

Let $T$ be a tilting $\Lambda$ module and $B={\rm End}_\Lambda\ T$.
Then $(\mathscr{T}(T),\mathscr{F}(T))$ is the torsion pair in
$\Lambda$-mod generated by $T$, where $\mathscr{T}(T)= T^{\bot}={\rm
gen}\ T$ and $\mathscr{F}(T)=\{X\in A-{\rm mod}\ |\ \ {\rm
Hom}_\Lambda(T,X)=0 \}$, the corresponding torsion pair in $B$-mod
is $(\mathscr{X}(T),\mathscr{Y}(T))$, where $\mathscr{X}(T)=\{X\in
B-{\rm mod}\ |\ \ T\otimes_B X=0 \}$ and $\mathscr{Y}(T)=\{Y\in
B-{\rm mod}\ |\ \ {\rm Tor}^B_1(T,Y)=0 \}$.

\vskip 0.2in

{\bf Lemma 2.1.} Take the notations as above. Then

\vskip 0.1in

(1)\ ${\rm Hom}_\Lambda(T,-): \mathscr{T}(T)\longrightarrow
\mathscr{Y}(T)$ is an equivalence functor;

\vskip 0.1in

(2)\ Let $M\in \mathscr{T}(T)$. Then ${\rm pd}_B\ {\rm
Hom}_\Lambda(T,M)\leq  {\rm pd}_\Lambda M $.     $\hfill\Box$

\vskip 0.2in

 Let ${\cal T}_{\Lambda}$ be the set of
all basic tilting $\Lambda$ modules up to isomorphism. According to
[10], we define the tilting quiver $\mathscr{K}(\Lambda)$ of
$\Lambda$ as the following. The vertices of $\mathscr{K}(\Lambda)$
are the elements of ${\cal T}_{\Lambda}$. There is an arrow
$T'\rightarrow T$ in $\mathscr{K}(\Lambda)$ if and only if
$T'=M\oplus X$ and $T=M\oplus Y$ with $X$ and $Y$ indecomposable
such that there is a short exact sequence $0\rightarrow X\rightarrow
E\rightarrow Y\rightarrow 0$ with $E\in{\rm add}\ M$.

\vskip 0.2in

Let $\mathcal{C}$ be a full subcategory of $\Lambda$-mod,
$C_{M}\in\mathcal{C}$ and $\varphi :C_M\longrightarrow M$ with
$M\in$ $A$-mod. Recall from [3], the morphism $\varphi$ is a right
$\mathcal{C}$-approximation of $M$ if the induced  morphism ${\rm
Hom}_{\Lambda}(C,C_{M})\longrightarrow {\rm Hom}_{\Lambda}(C,M)$ is
surjective for any $C\in\mathcal{C}$. A minimal right
$\mathcal{C}$-approximation of $M$ is a right
$\mathcal{C}$-approximation which is also a right minimal morphism,
i.e., its restriction to any nonzero summand is nonzero. The
subcategory $\mathcal{C}$ is called contravariantly finite if any
module $M\in$ $A$-mod admits a (minimal) right
$\mathcal{C}$-approximation. The notions of (minimal) left
$\mathcal{C}$-approximation and of covariantly finite subcategory
can be defined dually. It is well known that ${\rm add}\ M$ is both
a contravariantly finite subcategory and a covariantly finite
subcategory.

\vskip 0.2in

Let $M, N $ be two indecomposable ${\Lambda}$-modules. A path from
$M$ to $N$ in $\Lambda$-ind is a sequence of non-zero morphisms
$M=M_0\stackrel{f_1} \longrightarrow M_1\stackrel{f_2}
\longrightarrow\cdots \stackrel{f_t} \longrightarrow M_t =N$ with
all $M_i$ in ${\Lambda}$-ind. Following [13], we denote by $M\leq N$
the existence of such a path, and we say that $M$ is a predecessor
of $N$ (or that $N$ is a successor of $M$).

\vskip 0.2in

From now on, let $A=kQ$ be a finite dimensional hereditary algebra
over a field $k$, and let $Q_0=\{ 1,\cdots, n\}$ be the vertexes set
of $Q$. Recall from [1],  $A^{(1)}= ( \begin{array}{lr} A & 0\\
DA & A \end{array} )$ is said to be the duplicated algebra of $A$.
We know that $A^{(1)}$ contains two copies of $A$ given by $eA^{(1)}e$
and by $e'A^{(1)}e'$ respectively, where $e=( \begin{array}{lr} 1 & 0\\
0 & 0 \end{array} )$,  and $e'=( \begin{array}{lr} 0 & 0\\
0 & 1 \end{array} )$.  We denote the first one by $A$ and the second
one by $A'$. Then we denote by $Q'$ the quiver of $A'$, by $i'$ the
vertex of $Q'_0$ corresponding to $i\in Q_0$,  and by $e'_i$ the
corresponding idempotent.  Let ${S}_x, {P}_x, {I}_x$ denote
respectively the corresponding simple, indecomposable projective and
indecomposable injective module in $A^{(1)}$ corresponding to
 $x\in (Q_0\bigcup Q'_0)$. Note that $A$-ind can be embedded in
 $A^{(1)}$-ind, and $P_{x'}$ is an indecomposable projective-injective
 $A^{(1)}$ module for every $x'\in Q'_0$.

\vskip 0.2in

We denote by $\Sigma_0$ the set of all non-isomorphic indecomposable
projective $A$-modules and by $\Sigma_i$ the set of
$\Omega_{A^{(1)}}^{-i}\Sigma_0$. Note that $2 \leq {\rm gl.dim }\
A^{(1)}\leq 3$. Moreover, if $A$ is representation-infinite, then
${\rm gl.dim }\ A^{(1)}= 3$. (See [12])

\vskip 0.2in

 Let ${\cal L}_{A^{(1)}}$ be the left part of mod $A^{(1)}$.
 By definitions, ${\cal L}_{A^{(1)}}$ is the full subcategory
 of mod ${A^{(1)}}$ consisting of all indecomposable $A^{(1)}$-modules
 such that if $L$ is a predecessor of $M$, then the projective dimension ${\rm pd}\ L$
 of $L$ is at most one.

\vskip 0.2in

The following result is proved in [15] and will be used in our
further research.

\vskip 0.2in

{\bf Lemma 2.2.} Let $A=kQ$ be a finite dimensional hereditary
algebra over a field $k$ and $A^{(1)}$ be the duplicated algebra of
$A$. Then the tilting quiver $\mathscr{K}(A^{(1)})$ is connected.
$\hfill\Box$

\vskip 0.2in

{\bf Remark.}\  We should mention that $A^{(1)}$ was generalized to
$m$-replicated algebra $A^{(m)}$ far any integer $m\geq 1$ in [2],
and this kinds of algebras has been proved having closely
relationship with $m$-cluster categories,  and was extensively
investigated in [11, 12, 15, 16].

\vskip 0.2in

Throughout this paper, we follow the standard terminology and
notations used in the representation theory of algebras as in [4,
13].

\vskip 0.2in

\section {Global dimension of endomorphism algebras of tilting modules}

\vskip 0.2in

Let $A=kQ$ be a finite dimensional hereditary algebra over a field
$k$ and $A^{(1)}$ be the duplicated algebra of $A$. For convenience,
we denote by $\bar{P}_{A^{(1)}}$ the direct sum of indecomposable
projective-injective $A^{(1)}$ modules. In this section, we prove
that global dimension of the endomorphism algebras of tilting
$A^{(1)}$ modules is at most $3$.

\vskip 0.2in

Let $T\oplus \bar{P}_{A^{(1)}}$ be a basic tilting $A^{(1)}$ module
and $B=\rm{End}(T\oplus\bar{P}_{A^{(1)}})$. We know that $T\in {\rm
add}\ {\cal L}_{A^{(1)}}$ and $\delta(T)=\delta(A)=n$.

\vskip 0.2in

By [1] we know that $\rm{gl.dim}\ A^{(1)}\leq 3$ and $\rm{gl.dim}\
A^{(1)}= 3$ if $A$ is representation infinite. It is well known that
$\rm{gl.dim}\ A^{(1)}-1\leqslant \rm{gl.dim}\ B\leqslant
\rm{gl.dim}\ A^{(1)}+1$ which implies that $\rm{gl.dim}\ B\leq 4$.
However, we can prove the following surprising result.

\vskip 0.2in

{\bf Theorem 3.1.}\ {\it Take the notations as above. Then
$\rm{gl.dim}\ B\leq$ 3.}

\vskip 0.1in

{\bf Proof} \ Let $T= \bigoplus\limits_{i=1}^{n} T_i$ and
$\bar{P}_{A^{(1)}}= \bigoplus\limits_{i=1}^{n} \bar{P}_{i'}$. Let
$S$ be a simple $B$ module.

\vskip 0.1in

{\bf Case 1.} \ Assume that $S$ is the top of
$\rm{Hom}_{A^{(1)}}(T\oplus\bar{P}_{A^{(1)}},T_{i})$. Then we have
an exact sequence $0\rightarrow Y\rightarrow
\rm{Hom}_{A^{(1)}}(T\oplus \bar{P}_{A^{(1)}},T_{i})\rightarrow S
\rightarrow 0$. Note that $Y$ lies in $\mathscr{Y}(T\oplus
\bar{P}_{A^{(1)}})$ since $\mathscr{Y}(T\oplus \bar{P}_{A^{(1)}})$
is a torsion free class and
$\rm{Hom}_{A^{(1)}}(T\oplus\bar{P}_{A^{(1)}},T_{i})$ lies in
$\mathscr{Y}(T\oplus \bar{P}_{A^{(1)}})$. According to Lemma 1.1,
there exists $M\in \mathscr{T}(T\oplus \bar{P}_{A^{(1)}})$ such that
$Y=\rm{Hom}_{A^{(1)}}(T\oplus\bar{P}_{A^{(1)}},M)$, hence $M$ is a
predecessor of $T_i$ and $M\in\mathcal {L}_{A^{(1)}}$ since $T_i$
lies in $\mathcal {L}_{A^{(1)}}$. Therefor $\rm{pd}_{A^{(1)}}\ M\leq
1$ and by Lemma 2.1 again, we know that $\rm{pd}_B\ Y=\rm{pd}_B\
\rm{Hom}_{A^{(1)}}(T\oplus\bar{P} _{A^{(1)}},M)\leq
\rm{pd}_{A^{(1)}}\ M\leq 1$, which implies that $\rm{pd}_B\ S\leq
2$.

\vskip 0.1in

{\bf Case 2.} \  Let $S$ be the top of $\rm{Hom}_{A^{(1)}}(T\oplus
\bar{P}_{A^{(1)}},\bar{P}_{i^{'}})$. Then we have an exact sequence
$0\rightarrow Y\rightarrow \rm{Hom}_{A^{(1)}}(T\oplus
\bar{P}_{A^{(1)}},\bar{P}_{i^{'}})\rightarrow S \rightarrow 0$. By
using the same argument as in Case 1, we know that
$Y=\rm{Hom}_{A^{(1)}}(T\oplus\bar{P}_{A^{(1)}},M)$ with $M\in
\mathscr{T}(T\oplus \bar{P}_{A^{(1)}})$ such that $M$ is a
predecessor of $\Sigma_{2}$,  hence $\rm{pd}_{A^{(1)}}\ M\leq 2$.
According to Lemma 2.1, we know that ${\rm pd}_B\ Y={\rm pd}_B\ {\rm
Hom}_{A^{(1)}}(T\oplus\bar{P} _{A^{(1)}},M)\leq {\rm pd}_{A^{(1)}}\
M\leq 2$, which implies that $\rm{pd}_B\ S\leq 3$. This proves that
$\rm{gl.dim}\ B\leq 3$.         $\hfill\Box$

\vskip 0.2in

\section {Embedding of the tilting quiver}

\vskip 0.2in

Let $A=kQ$ be a finite dimensional hereditary algebra over a field
$k$ and $A^{(1)}$ be the duplicated algebra of $A$. In this section,
we investigate the relationship between the tilting quivers of $A$
and of $A^{(1)}$.

\vskip 0.2in

{\bf Theorem 4.1.}\ {\it Let $\mathscr{K}(A)$ (resp.
$\mathscr{K}(A^{(1)})$) be the tilting quiver of $A$ (resp.
$A^{(1)}$). Then there is an arrow $T^{'}\rightarrow T$ in
 $\mathscr{K}(A)$ if and only if $T^{'}\oplus
\bar{P}_{A^{(1)}}\rightarrow T\oplus \bar{P}_{A^{(1)}}$ is an arrow
in  $\mathscr{K}(A^{(1)})$.}

\vskip 0.1in

{\bf Proof}\ \  Let $T$ is a tilting $A$ module. It is easy to see
that $T\oplus \bar{P}_{A^{'}}$ is a tilting $A^{(1)}$ module.

Assume that $T^{'}\rightarrow T$ is an arrow in $\mathscr{K}(A)$,
then there is an almost tilting $A$ module $M$ such that
$T^{'}=M\oplus X$ and $T=M\oplus Y$ with $X$ and $Y$ are
indecomposable. Moreover,  there is an exact sequence $0\rightarrow
X\stackrel{f}{\rightarrow} E\stackrel{g}{\rightarrow} Y\rightarrow
0$ is an exact sequence  with $E\in\rm{add}\ M$, such that $f$ is a
left minimal $\rm{add}\ M$-approximation and that $g$ is a right
minimal $\rm{add}\ M$-approximation.

It follows that $T^{'}\oplus\bar{P}_{A^{(1)}}=M\oplus X\oplus
\bar{P}_{A^{(1)}}$ and $T\oplus \bar{P}_{A^{(1)}}=M\oplus
Y\oplus\bar{P}_{A^{'}}$ are tilting $A^{(1)}$ modules, and $g$, $f$
are also minimal ${\rm add}\ M\oplus\bar{
P}_{A^{(1)}}$-approximation, since $ \bar{P}_{A^{(1)}}$ is a
projective-injective module. Hence $T^{'}\oplus
\bar{P}_{A^{(1)}}\rightarrow T\oplus \bar{P}_{A^{(1)}}$ is an arrow
in $\mathscr{K}(A^{(1)})$.

The converse can be proved similarly. This completes the proof.
$\hfill\Box$

\vskip 0.2in

{\bf Theorem 4.2.}\ {\it  Each point in $\mathscr{K}(A^{(1)})$ has n
arrows connected.}

\vskip 0.1in

{\bf Proof}\ \ Let $T\oplus\bar{P}_{A^{(1)}}$ be a basic tilting
$A^{(1)}$ module. Then $\delta(T)=\delta(P_{A^{(1)}})=n$. Assume
that $T=\bigoplus\limits^{n}_{i=1} T_{i}$, and let
$T[i]=\bigoplus\limits_{j\neq i} T_{j}$. Then
$T[i]\oplus\bar{P}_{A^{(1)}}$ is an almost tilting $A^{(1)}$ module.

According to [14], we know that $T[i]\oplus\bar{P}_{A^{(1)}}$ has
exactly two non-isomorphic complements with projective dimension at
most 1, and one of them is $T_{i}$.

Note that if $T_{i}$ is the source complement, then there exists an
arrow $T\oplus \bar{P}_{A^{(1)}}\rightarrow *$. Otherwise, there is
an arrow $*\rightarrow T\oplus P_{A^{(1)}}$. This implies that there
are exactly  n arrows  connected with $T\oplus \bar{P}_{A^{(1)}}$.
The proof is completed.            $\hfill\Box$

\vskip 0.2in

{\bf Theorem 4.3.}\ {\it Let $M$ be a basic almost tilting $A$
module. Then $({\bf dim}\ M){i}=0$ if and only if
$M\oplus\tau_{A^{(1)}}^{-1}I_{i}\oplus \bar{P}_{A^{(1)}}$ is a
tilting $A^{(1)}$ module.}

\vskip 0.1in

{\bf Proof}\ \ Note that ${\rm pd}_{A^{(1)}}\ M\leq 1$ and
$\tau_{A^{(1)}}^{-1}I_{i}\in\mathcal {L}_{A^{(1)}}$, it follows that
$$
\rm{pd}_{A^{(1)}}\ (M\oplus\tau_{A^{(1)}}^{-1}I_{i}
\oplus\bar{P}_{A^{(1)}})\leq 1.
$$
We have that
$${\rm Ext}^{1}_{A^{(1)}}(M,\tau_{A^{(1)}}^{-1}I_{i})\cong D {\rm
Hom}_{A^{(1)}}(\tau_{A^{(1)}}^{-1}I_{i},\tau_{A^{(1)}} M)=0,
$$
and that
$$
{\rm Ext}^{1}_{A^{(1)}}(\tau_{A^{(1)}}^{-1}I_{i},M)\cong D{\rm
Hom}_{A^{(1)}}(M,I_{i})=0,
$$
hence
$${\rm Ext}^{1}_{A^{(1)}}(M\oplus\tau_{A^{(1)}}^{-1}I_{i}\oplus\bar{P}_{A^{1)}},
M\oplus\tau_{A^{(1)}}^{-1}I_{i}\oplus\bar{P }_{A^{1)}})= 0,
$$
then $M\oplus\tau_{A^{(1)}}^{-1}I_{i}\oplus\bar{P}_{A^{1)}}$ is a
tilting $A^{(1)}$ module, since
$$
\delta(M\oplus\tau_{A^{(1)}}^{-1}I_{i}\oplus\bar{P}_{A^{'}})=\delta(A^{(1)})=2n.
$$

\vskip 0.1in

Conversely, if $M\oplus\tau_{A^{(1)}}^{-1}I_{i}\oplus
\bar{P}_{A^{(1)}}$ is a tilting $A^{(1)}$ module, then
$$
{\rm Ext}^{1}_{A^{(1)}}(\tau_{A^{(1)}}^{-1}I_i,M)\cong D{\rm
Hom}_{A^{(1)}}(M,I_i)=0,
$$
this implies that ${\rm Hom}_{A}(M,I_i)=0$  and $({\bf dim}\
M)_{i}=0$. The proof is completed. $\hfill\Box$

\vskip 0.2in

The following corollary can be proved easily.

\vskip 0.2in

{\bf Corollary 4.4.}\ {\it Let $M$ be an almost tilting $A$ module
and $M$ is not sincere, then the dimension vector ${\bf dim}\ M$ of
$M$ has exactly one component equals to 0.}

\vskip 0.1in

{\bf Proof}\  Assume by contrary that there are two or more
different components of ${\bf dim}\ M$ equal to zero. That is, there
are $i\neq j$ such that $({\bf dim}\ M)_i=({\bf dim}\ M)_j=0$. By
using the method in the proof of Theorem 4.3, we know that
$M\oplus\tau_{A^{(1)}}^{-1}I_{i}\oplus\tau_{A^{(1)}}^{-1}I_{j}\oplus\bar{
P}_{A^{'}}$ is a tilting $A^{(1)}$ module, then
$\delta(M\oplus\tau_{A^{(1)}}^{-1}I_{i}\oplus\tau_{A^{(1)}}^{-1}I_{j}
\oplus\bar{P}_{A^{'}})=2n+1$, which is a contradiction. $\hfill\Box$

\vskip 0.2in

\section {Applications of the embedding theorem}

\vskip 0.2in

Let $A=kQ$ be a finite dimensional hereditary algebra over a field
$k$ and $A^{(1)}$ be the duplicated algebra of $A$. In this section,
we give new proofs for some results of D.Happel and L.Unger by using
embedding theorem, and obtain a partially positive answer to the
conjecture of D.Happel and L.Unger in [10], which says that every
connected component in $\mathscr{K}({A})$ has finite non-saturated
points. We also prove that the number of arrows in
$\mathscr{K}({A})$ is a constant which does not depend on the
orientation of $Q$ if $Q$ is Dynkin type.

\vskip 0.2in

The following proposition is the main result of [7], we give a new
proof by using embedding theorem in Section 3.

\vskip 0.2in

{\bf Proposition 5.1.}$^{[7]}$ \ {\it Let $M$ be an almost tilting
$A$ module. If $M$ is sincere, then $M$ has two non-isomorphic
indecomposable complements, and if $M$ is non-sincere, then $M$ has
exactly one complement.}

\vskip 0.1in

{\bf Proof}\ \  Note that $M\oplus\bar{P}_{A^{(1)}}$ is an almost
tilting $A^{(1)}$ module,  according to [15], we know that
$M\oplus\bar{P}_{A^{(1)}}$ has two non-isomorphic indecomposable
complements $X, Y\in {\rm ind}\ A\cup\{\tau_{A^{(1)}}^{-1}I_{x}|x\in
Q\}$.

\vskip 0.1in

If $M$ is sincere, then $X, Y\in\rm{ind}A$. Otherwise, we may assume
that $X=\tau_{A^{(1)}}^{-1}I_{i}$, according to Theorem 4.3, $({\bf
dim}\ M)_{i}=0$ which means $M$ not sincere.

\vskip 0.1in

If $M$ is non-sincere, by Corollary 4.4 we know that
$M\oplus\bar{P}_{A^{(1)}}$ has exactly one complement looking like
$\tau_{A^{(1)}}^{-1}I_{i}$, and the other complement must lie in
$\rm{ind}A$ which is also the only complement for the almost tilting
$A$ modules $M$.      $\hfill\Box$

\vskip 0.2in

Recall from [10], let $T\in\mathscr{K}(A)$. We denote by $s(T)$
(resp. e(T)) the number of arrows starting (resp. ending) at $T$ in
$\mathscr{K}(A)$, then $\sigma(T) = s(T)+e(T)\leq \delta (A)=n$. We
say that $T$ is saturated if $\sigma(T) = n$. The following result
is stated as Proposition 3.2 in [10], and we provide a new proof
here.

\vskip 0.2in

{\bf Proposition 5.2.}$^{[10]}$ \ {\it Let $T$ be a basic tilting
$A$ module, then the point $T$ in the tilting quiver
$\mathscr{K}(A)$ of $A$ is saturated if and only if $({\bf dim}\
T)_{i}\geq 2, \forall i\in Q_0$.}

\vskip 0.1in

{\bf Proof}\ \ Assume that $T$ is saturated and there is some $i\in
Q_0$ with $({\bf dim}\ T)_{i}= 1$, then there must be an
indecomposable summand $T_{k}$ of $T$ such that $({\bf dim}\
T_{k})_{i} = 1$. So $T[k]$ is non-sincere since the $i^{th}$
component of ${\bf dim}\ T[k]$ is 0. According to Proposition 4.1,
there is only one complement for $T[k]$ in $A$-mod. This means that
$T$ is not saturated, and we get a contradiction.

\vskip 0.1in

Conversely,  If $({\bf dim}\ T)_{i}\geq 2$  for all $i\in Q_0$ and
$T$ is not saturated, then we know that there exists at least one
$T[k]$, in ${\rm mod}\ A$, which has the unique complement $T_k$,
hence $T[k]$ is non-sincere. We may assume that $({\bf dim}\
T[k])_{i}=0$, according to Theorem 4.3, $T[k]\oplus
\bar{P}_{A^{(1)}}$ has a complement $\tau_{A^{(1)}}^{-1}I_{j}$ in
${\rm mod}\ A^{(1)}$. It follows that $T[k]\oplus \bar{P}_{A^{(1)}}$
has two complements $X=T_{k}$ and $\tau_{A^{(1)}}^{-1}I_{j}$, which
means that there is an exact sequence $0\rightarrow X\rightarrow
E\rightarrow \tau_{A^{(1)}}^{-1}I_{j}\rightarrow 0$ with $E\in
\rm{add}(T[k]\oplus \bar{P}_{A^{(1)}})$. Applying ${\rm
Hom}_{A^{(1)}}(-,I_{j})$ we obtain the following exact sequence
$$
{\rm Hom}_{A^{(1)}}(E,I_{i}) \rightarrow {\rm
Hom}_{A^{(1)}}(X,I_{j}) \rightarrow {\rm
Ext}_{A^{(1)}}^{1}(\tau_{A^{(1)}}^{-1}I_{j},I_{j}) \rightarrow 0.
$$
${\rm Hom}_{A^{(1)}} (E,I_{j})=0$ since ${\rm
Hom}_{A^{(1)}}(T[k],I_{j})=0$ and ${\rm Hom}_{A^{(1)}}
(\bar{P}_{A^{(1)}},I_{j})=0$, hence
$$
({\bf dim}\ X)_{j}= {\rm dim\ Hom}_{A^{(1)}}(X,I_{j})={\rm dim\
Ext}_{A^{(1)}}^{1}(\tau_{A^{(1)}}^{-1}I_{j},I_{j})= {\rm dim}\ D
{\rm Hom}(I_{j},I_{j})=1.
$$
It follows that $({\bf dim}\ T)_j=({\bf dim}\ T[k])_j+({\bf dim}\
X)_j=1$, which contradicts with the assumption. $\hfill\Box$

\vskip 0.2in

{\bf Corollary 5.3.}$^{[10]}$ \ {\it Let $A=kQ$ be a finite
dimensional hereditary algebra over a field $k$. Then $A$ and $DA$
are not saturated in the tilting quiver $\mathscr{K}(A)$.}

\vskip 0.1in

{\bf Proof} \ Let $i$ be a source vertex of $Q_0$. Then $({\bf dim}\
\bigoplus\limits_{{j\in Q_0},\  {j\neq i}} P_{j})_{i}=0$, hence $A$
is not saturated. That $DA$ is not saturated can be proved dually.
$\hfill\Box$

\vskip 0.2in

We give a very different proof for Theorem 3.5 in [10] as following.

\vskip 0.2in

{\bf Proposition 5.4.}$^{[10]}$ \ {\it  Let $A=kQ$ be a finite
dimensional hereditary algebra over a field $k$. Then each connected
component in the tilting quiver $\mathscr{K}(A)$ has a non-saturated
point.}

\vskip 0.1in

{\bf Proof}\ \ If $\mathscr{K}({A})$ is connected, it is easy to see
that $A$ is one of non-saturated point in $\mathscr{K}({A})$. Now,
we assume that $\mathscr{K}({A})$ is not connected. If
$\mathscr{K}({A})$ has one component such that every point is
saturated, according to Proposition 5.1 and Proposition 4.2,
$\mathscr{K}(A)$ can be embedded into $\mathscr{K}(A^{(1)})$ and the
only change is that every basic tilting $A$ module $T$ is replaced
by $T\oplus\bar{P}_{A^{(1)}}$ and the arrows keep no changes. This
implies that the component, which every point is saturated, is
isolated. In particular, $\mathscr{K}(A^{(1)})$ has at least two
components, which is contradict with Lemma 2.2. This completes the
proof.            $\hfill\Box$

\vskip 0.2in

Let $A=kQ$ be a finite dimensional hereditary algebra over a field
$k$. D.Happel and L.Unger in [10] conjectured that every connected
component in $\mathscr{K}({A})$ has finite non-saturated points. The
following theorem gives a partially positive answer to this
conjecture.

\vskip 0.2in

{\bf Theorem 5.5.} \ {\it  Let $A=kQ$ be a finite dimensional
hereditary algebra over a field $k$. If $Q$ is either Dynkin  or
Euclidean type, then every connected component of $\mathscr{K}({A})$
has finite non-saturated points.}

\vskip 0.2in

{\bf Proof}\ \  If $A$ is a $\rm{Dynkin}$ type, then
$\mathscr{K}(A)$ is a finite quiver and our consequence is true. Now
we assume that $Q$ is an Euclidean type. Let $T$ be a  non-saturated
point in $\mathscr{K}(A)$. Then ${\bf dim}\ T$ has at least one
component equal to 1. We denote by $\Delta$ the set of non-saturated
points in $\mathscr{K}(A)$, and we divide  $\Delta$ into different
parts and put $\Delta_i=\{ T\in \Delta\ |\ ({\bf dim}\ T)_i=1\}$.

We claim that $\Delta_i$ is a finite set for $1\leq i\leq n$. In
fact, $\forall T\in\Delta_i$, we know that $({\bf dim}\ T)_i=1$. Let
$T=\bigoplus\limits_{i=1}^n \ T_i$. Then there is a $T_k$ with
$({\bf dim}\ T_k)_i=1$. Let $Q_{(i)}$ be the quiver by removing the
vertex $i$ from $Q_0$ and removing all the arrows connected with
$i$. Then $T[k]$ can be regarded as a basic tilting $kQ_{(i)}$
module, and $kQ_{(i)}$ is representation-finite, hence $\Delta_i$ is
a finite set, it follows that $\Delta=\bigcup\limits_{i=1}^n
\Delta_i$ is also a finite set. This completes the proof.
$\hfill\Box$

\vskip 0.2in

{\bf Theorem 5.6.} \ {\it Let $A=kQ$ and $Q$ be $\rm{Dynkin}$ type.
Then the number of arrows in $\mathscr{K}(A)$ is a constant and does
not depend on the orientation of Q.}

\vskip 0.1in

{\bf Proof}\ \ According to Theorem 4.1,  we know that
$\mathscr{K}(A)$ can be embedded into $\mathscr{K}(A^{(1)})$. We
denote by $\widehat{\mathscr{K}(A)}$ the full subquiver of
$\mathscr{K}(A)$ in $\mathscr{K}(A^{(1)})$. Note that
$\widehat{\mathscr{K}(A)}$ has the same vertices as
$\mathscr{K}(A)$, and every vertex in $\widehat{\mathscr{K}(A)}$
connected with $n$ arrows. Let $s$ be the number of basic tilting
$A$ modules, it is well known that $s$ is a fixed number which is
independent of the orientation of $Q$.

Let $Q_{\widehat{x}}=Q\setminus \{x\}$ be the quiver obtained from
$Q$ with a vertex $x\in Q_0$ removed. Then $kQ_{\widehat{x}}$ is a
representation finite hereditary algebra. We denote by $m_x$ the
number of basic tilting $kQ_{\widehat{x}}$ modules, then $m_x$ is a
fixed number which does not depend on the orientation of $Q$.

Let $m$ be the number of arrows should be added in order to get
$\widehat{\mathscr{K}(A)}$ from $\mathscr{K}(A)$. Note that every
tilting $kQ_{\widehat{x}}$ module is a non-sincere almost tilting
$A$ module, and there is one corresponding arrow in
$\widehat{\mathscr{K}(A)}\backslash \mathscr{K}(A)$. On the other
hand, every arrow in $\widehat{\mathscr{K}(A)}\backslash
\mathscr{K}(A)$ corresponding to one almost tilting $A$ module which
can be seen as tilting $kQ_{\widehat{x}}$ module for some $x\in
Q_0$. According to Theorem 4.3, $m= \sum\limits_{x\in Q_0} m_x$ is a
fixed number.

Let $t$ be the number of arrows in $\mathscr{K}(A)$. According to
Theorem 4.2, we have an equation $2t+m=ns$, hence $t=\frac{ns-m}{2}$
is a fixed number which does not depend on the orientation of $Q$,
that is, $t$ is a constant. The proof is finished. $\hfill\Box$

\vskip 0.2in

{\bf Remark.} \ \ Theorem 5.6 is more general than the result in
[14] which stands only for $A_n$ and $D_n$ type, and our proof is
different and simpler.

\vskip 0.2in

{\bf Acknowledgements.}\  This paper was finally written while the
second author visited the Department of Mathematics of Hong Kong
University of Science and Technology. He would like to thank
Prof.Jing-Song Huang for his helpful suggestions and kind
hospitality, and he also would like to thank the Department of
Mathematics for facilities.

\begin{description}

\item{[1]}\ I.Assem, T.Br$\ddot{{\rm u}}$stle, R.Schiffler, G.Todorov,
   Cluster categories and duplicated algebras.
   {\it J. Algebra }, 305(2006), 548-561.

\item{[2]}\ I.Assem, T.Br$\ddot{\rm u}$stle, R.Schiffer, G.Todorov,
 $m$-cluster categories and $m$-replicated algebras.  {\it J.
Pure and Appl. algebra},  212(2008), 884-901.

\item{[3]}\ M.Auslander, I.Reiten,  Applications of contravariantly finite
subcategories.  {\it Adv. Math.}, 86(1991), 111-152.

\item{[4]}\ M.Auslander, I.Reiten, S.O.Smal$\phi$,  Representation
Theory of Artin Algebras.  Cambridge Univ. Press, 1995.

\item{[5]}\ A. Buan, R.Marsh, M.Reineke, I.Reiten, G.Todorov,
Tilting theory and cluster combinatorices, {\it  Adv. Math.,} 204(2)(2006), 572-618.

\item{[6]}\ D.Happel, C.M.Ringel,  Tilted algebras. {\it Trans.Amer.Math.Soc.,}
274(1982), 399-443.

\item{[7]}\  D.Happel, L.Unger, Almost complete tilting modules.
 {\it Proc.Amer.Math. Soc.}, 107(1989), 603-610.

\item{[8]}\ D.Happel, L.Unger,  Partial tilting modules and covariantly finite
subcategories. {\it Comm.Algebra},  22(1994), 1723-1727.

\item{[9]}\ D.Happel, L.Unger,  Reconstruction of path algebras from their posets of tilting modules.
{\it Trans.Amer.Math.Soc.,} 361(2009), 3633-3660.

\item{[10]}\ D.Happel, L.Unger, On the quiver of tilting modules. J.Algebra, 284(2005), 857-868.

\item{[11]}\ X.Lei, H.Lv, S.Zhang, Complements to the almost complete tilting
$A^{(m)}$-modules.  {\it Comm.Algebra},  37(2009), 1719-1728.

\item{[12]}\ H.Lv, S.Zhang, Global dimensions of endomorphism algebras for
generator-cogenerators over m-replicated algebras.  {\it
Comm.Algebra},  39(2011), 560¨C571.

\item{[13]}\ C.M.Ringel,  Tame algebras and integral quadratic forms.
{\it Lecture Notes in Math. 1099.}  Springer-Verlag, Berlin,
Heidelberg, New York, 1984.

\item{[14]}\ K.Ryoichi,  The number of arrows in the quiver of
tilting modules over a path algebra of type A and D. Preprint 2011,
arXiv:1101.4747.

\item{[15]}\ S.Zhang,  Tilting mutation and duplicated
algebras.  {\it Comm.Algebra}, 37(2009), 3516-3524.

\item{[16]}\ S.Zhang, Partial tilting modules over m-replicated
algebras. {\it J.Algebra}, 323 (2010), 2538-2546

\end{description}

\end{document}